\documentclass[english]{article}
\usepackage[T1]{fontenc}
\usepackage[latin9]{inputenc}
\usepackage{amsthm}
\usepackage{amsmath}

\makeatletter
\numberwithin{equation}{section}
\numberwithin{figure}{section}
\theoremstyle{plain}
\newtheorem{thm}{\protect\theoremname}
  \theoremstyle{remark}
  \newtheorem{rem}[thm]{\protect\remarkname}
  \theoremstyle{plain}
  \newtheorem{lem}[thm]{\protect\lemmaname}
  \theoremstyle{plain}
  \newtheorem{cor}[thm]{\protect\corollaryname}

\makeatother

\usepackage{babel}
  \providecommand{\corollaryname}{Corollary}
  \providecommand{\lemmaname}{Lemma}
  \providecommand{\remarkname}{Remark}
\providecommand{\theoremname}{Theorem}

\begin{document}

\title{Characterizing the second smallest eigenvalue of the normalized Laplacian
of a tree}

\author{Israel Rocha}
\maketitle
\begin{abstract}
In this paper we show a monotonicity theorem for the harmonic eigenfunction
of $\lambda_{1}$ of the normalized Laplacian over the points of articulation
of a graph. We introduce the definition of Perron component for the
normalized Laplacian matrix of a graph and show how its second smallest
eigenvalue can be characterized using this definition.
\end{abstract}

\section{Main Concepts}

As usual, in this paper a graph is a pair of sets $G=(V,E)$, where
the elements of $E$ are subsets of two elements of $V$. The elements
of $V$ are vertices of the graph and the elements of $E$ are its
edges.

Given a graph $G=(V,E)$ on $n$ vertices, the \emph{normalized Laplacian
matrix }of $G$ is the matrix of order $n$ $\mathcal{L}(G)$ given
by 
\[
\mathcal{L}(v_{i},v_{j})=\begin{cases}
1, & \textrm{}v_{i}=v_{j}\textrm{ and }d_{v_{i}}\neq0;\\
-\frac{1}{\sqrt{d_{v_{i}}d_{v_{j}}}}, & \textrm{whenever }v_{i}\textrm{ and }v_{j}\textrm{ are adjacent;}\\
0,\textrm{ } & otherwise.
\end{cases}
\]
Also the \textit{Laplacian matrix} of $G$ is the matrix of order
$n$ given by 

\[
L(v_{i},v_{j})=\begin{cases}
d(v_{i}), & \textrm{}v_{i}=v_{j};\\
-1, & \textrm{whenever }v_{i}\textrm{ and }v_{j}\textrm{ are adjacent;}\\
0,\textrm{ } & otherwise.
\end{cases}
\]

In the survey \cite{key-6}, some known results about on Laplacian
matrix are exhibit. Fiedler in \cite{key-4}, has shown that a graph
is connected if and only if the second smallest Laplacian eigenvalue
is positive. This eigenvalue is called \textit{algebraic connectivity}
and plays a fundamental role in the field of Spectral Graph Theory. 

Throughout this paper, $G$ does not have isolated vertices. In that
case $D$ is invertible and $\mathcal{L}$ and $L$ are related by
the formula 
\[
\mathcal{L}=D^{-\frac{1}{2}}LD^{-\frac{1}{2}}.
\]

We consider the normalized laplacian matrix $\mathcal{L}$ of a tree
and, following the notation of \cite{key-1}, we denote the eigenvalues
of $\mathcal{L}$ by $0=\lambda_{0}\leq\lambda_{1}\leq\ldots\leq\lambda_{n-1}$.
Let $g$ denote an function which assigns to each vertex $v$ of $G$
a real value $g(v)$. We can view $g$ as a column vector and whenever
$\mathcal{L}g=\lambda g$ we call $g$ a eigenfunction of $\mathcal{L}$.
We define the harmonic eigenfunction of $\lambda$ as $f=D^{-\frac{1}{2}}g$
.

The following result, that we can find at \cite{key-2}, concerns
the harmonic eigenfunction of $\lambda_{1}$.
\begin{thm}
\label{thm:valuationHarmonicEigenfunction}Let $G$ be a connected
graph and $\mathcal{L}(G)$ be its normalized Laplacian matrix. Let
$f$ be a harmonic eigenfunction corresponding to $\lambda_{1}$ and
$v$ be a cut vertex of $G$, let $G_{0},G_{1},...,G_{r}$ be all
connected components of the graph $G\backslash v$. Then:

(1) If $f(v)>0$ then exactly one of the components $G_{i}$ contains
a vertex negatively valuated by $f.$ For all vertices $u$ in the
remaining components $f(u)>f(v).$

(2) If $f(v)=0$ and there exists a component $G_{i}$ containing
both positively and negatively valuated vertices, then there is exactly
one such component, all remaining components being zero valuated.

(3) If $f(v)=0$ and no component contains both positively and negatively
valuated vertices then each component contains either only positively
valuated, or negatively valuated, or zero valuated vertices.
\end{thm}
We notice that this result is similar to the result of Fiedler \cite{key-4}
where the eigenvector associated with the algebraic connectivity was
considered. In this paper we show a property of the harmonic eigenfunction
of $\lambda_{1}$ over the points of articulation which, likewise
in \cite{key-4}, enables us to classify every graph in two distinct
families. Also, we introduce the definition of Perron component for
the normalized Laplacian matrix of a graph, and using this we can
provide a characterization for $\lambda_{1}$ in therms of this definition.
Moreover, we introduce the notion of \emph{normalized bottleneck matrix}
of a branch of a tree which allow us to easily describe $\lambda_{1}$.
Furthermore, we shall perform a more careful analysis on the structure
of normalized bottleneck matrices in order to understand how $\lambda_{1}$
behaves when we change the structure of trees.

\section{Monotonicity Theorem}

In this section we show an interesting property of the harmonic eigenfunction
of $\lambda_{1}$ over the points of articulation of a graph. We shall
provide a monotonicity theorem for such harmonic eigenfunction. This
enable us to classify every graph in two distinct families.

First, a block of a graph is a maximal induced connected subgraph
not containing a point of articulation. A path is said to be pure
if it contains at most two points of articulation of each block.
\begin{thm}
\label{thm:monoHarmonic}Let $G$ be a connected graph and let $f$
be the harmonic eigenfunction for $\lambda_{1}$. Then only one of
the following cases can occur:

\textbf{Case 1 }There is no mixed block. In this case, there is a
unique point of articulation $z$ having $f(z)=0$ and a nonzero neighbor.
Each block (with the exception of the vertex $z$) is either a positive
block, or a negative block, or a zero block. Let $P$ be a pure path
which starts at $z$. Then the $f$ at the points of articulation
(with the exception of $z$) form either an increasing, or decreasing,
or a zero sequence. Every path containing both positive and negative
vertices passes through $z$.

\textbf{Case 2} There is a unique block $B_{0}$ which is mixed. In
this case, each remaining block is positive, negative or null. Moreover,
each pure path $P$ starting in $B_{0}$ and containing only one vertex
$v\in B_{0}$ has the property that $f$ at the points of articulation
contained in $P$ form either an increasing, or decreasing, or a zero
sequence according to whether $f(v)>0$, \textup{$f(v)<0$} or \textup{$f(v)=0$}.
In the last case $f\equiv0$ along the path.\end{thm}
\begin{proof}
First, for case 1, if no block is mixed, since $\sum d_{v}f(v)=0$,
there is a path containing both positive and negative vertices. We
claim that $P$ has a vertex $z$ with $f(z)=0$ and a nonzero neighbor.
Indeed, the intersection of blocks has only articulation points and
no block is mixed, it follows that exists such vertex. Thus, it follows
from Theorem \ref{thm:valuationHarmonicEigenfunction} that part (3)
must occurs. Therefore, there is no other vertex $v\neq z$ having
$f(v)=0$ and a nonzero neighbor. This shows the first part of case
1.

Now, if $P$ contains another vertex $v$ with $f(v)=0$, part (3)
of Theorem \ref{thm:valuationHarmonicEigenfunction} ensures that
$f\equiv0$ over the vertex of $P$. On the other hand, if $P$ has
a vertex $v$ with $f(v)\neq0$ then part (1) of Theorem \ref{thm:valuationHarmonicEigenfunction},
we obtain that $f$ does not change sign neither vanish over $P$.
Denote by $z=v_{0},v_{1},\ldots,v_{s}$ the points of articulation
at $P$ in the order they appear. If $f(v)>0$, then by part (1) of
Theorem \ref{thm:valuationHarmonicEigenfunction} we obtain $f(v_{i})<f(v_{i+1})$,
$i=0,\ldots,s-1$. If $f(v)<0$, then the same argument applied to
the eigenfunction $-f$, shows that this form a decreasing sequence. 

Now we proceed proving case 2. If $G$ has only one block, then we
are done. Otherwise, denote by $B_{1}$ some other block different
of $B_{0}$. In this case, there is a articulation point $v$ separating
them. Let $G_{0},G_{1},...,G_{r}$ be the connected components of
$G\backslash v$, where $G_{0}$ contains $B_{0}$ and $G_{1}$ contains
$B_{1}$. If $f(v)>0$ (or $f(v)<0$), then by part (1) of Theorem
\ref{thm:valuationHarmonicEigenfunction}, we obtain that $f$ has
the same sign over $G_{1}$. If $f(v)=0$, then using part (2) of
Theorem \ref{thm:valuationHarmonicEigenfunction}, we obtain that
$f\equiv0$ over $G_{1}$. This completes the first part of case 2.

Finally, denote by $v=v_{0},v_{1},\ldots,v_{s}$ the points of articulation
at $P$ in the order they appear. If $f(v)>0$, then by part (1) of
Theorem \ref{thm:valuationHarmonicEigenfunction} we obtain $f(v_{i})<f(v_{i+1})$,
$i=0,\ldots,s-1$. If $f(v)<0$, then the same argument applied to
the eigenfunction $-f$, shows that this form a decreasing sequence.
If $f(v)=0$, then using part (2) of Theorem \ref{thm:valuationHarmonicEigenfunction},
we obtain that $f\equiv0$ over the vertices of $P$. This concludes
the proof.\end{proof}
\begin{rem}
\label{remark1}Since $sign(f(v))=sign(g(v))$ for each vertex $v$
at $G$, we can provide the following result which is straightforward
from Theorem \ref{thm:valuationHarmonicEigenfunction}.\end{rem}
\begin{thm}
\label{thm:valuationEigenfunction}Let $G$ be a connected graph and
$\mathcal{L}(G)$ be its normalized Laplacian matrix. Let $g$ be
a eigenfunction corresponding to $\lambda_{1}$ and $v$ be a cut
vertex of $G$, let $G_{0},G_{1},...,G_{r}$ be all connected components
of the graph $G\backslash v$. Then:

(1) If $g(v)>0$ then exactly one of the components $G_{i}$ contains
a vertex negatively valuated by $g.$

(2) If $g(v)=0$ and there exists a component $G_{i}$ containing
both positively and negatively valuated vertices, then there is exactly
one such component, all remaining components being zero valuated.

(3) If $g(v)=0$ and no component contains both positively and negatively
valuated vertices then each component contains either only positively
valuated, or negatively valuated, or zero valuated vertices.
\end{thm}
We notice that in part (1) of Theorem \ref{thm:valuationEigenfunction}
unlike Theorem \ref{thm:valuationHarmonicEigenfunction}, we can not
guarantee that $g(u)>g(v)$, since de degree $d_{v}$ could be much
larger than $d_{u}$. Remark \ref{remark1} and Theorem \ref{thm:monoHarmonic}
give us the following result.
\begin{thm}
\label{thm:monoEigenfunction}Let $G$ be a connected graph and let
$g$ be the eigenfunction for $\lambda_{1}$. Then only one of the
following cases can occur:

\textbf{Case 1 }There is no mixed block. In this case, there is a
unique point of articulation $z$ having $g(z)=0$ and a nonzero neighbor.
Each block (with the exception of the vertex $z$) is either a positive
block, or a negative block, or a zero block. 

\textbf{Case 2} There is a unique block $B_{0}$ which is mixed. In
this case, each remaining block is positive, negative or null. 
\end{thm}
Henceforth, we use Theorem \ref{thm:monoEigenfunction} as it describes
directly the valuation of an eigenvector at the vertices of $G$.

\section{Characterizing the Second Smallest Eigenvalue}

Despite of giving classification of graphs and a good insight about
the behavior of the harmonic eigenfunction, Theorems \ref{thm:monoHarmonic}
and \ref{thm:monoEigenfunction} do not give us information about
$\lambda_{1}$ itself. However, we can provide a alternative characterization
for cases 1 and 2 such that information about $\lambda_{1}$ arises.

More precisely, in this section we are interested in describe $\lambda_{1}$
in terms of the Perron value of special matrices. This results were
inspired by \cite{key-3}. 

Consider the normalized Laplacian matrix $\mathcal{L}(G)$ for a graph
$G$. The relation between the matrix $L(G)$ and $\mathcal{L}(G)$
is well-known, and it is given by 
\[
\mathcal{L}=D^{-\frac{1}{2}}LD^{-\frac{1}{2}},
\]
where $D$ is the degree matrix. 

We shall denote by $M_{k}$, the principal submatrix of $M$ formed
by removing the $k-$th row and column of $M$. Now, consider the
matrix $\mathcal{L}_{k}$. We call normalized bottleneck matrix of
a component $C$ at $k$, the corresponding block at $\mathcal{L}_{k}^{-1}$. 

If we call $N=\mathcal{L}_{k}^{-1}(C)$ the normalized bottleneck
matrix of the component $C$, we say that it is a Perron component
if it has largest $\rho(N)$ among all components.

Let $T$ be a tree. We call branch of $T$ at $k$ some of the connected
components of $T-k$ obtained from $T$ by deleting the vertex $k$
and its edges. If $T$ satisfies case 1 of Theorem \ref{thm:monoEigenfunction}
then we say $T$ is a Type 1 tree. If $T$ satisfies case 2 of Theorem
\ref{thm:monoEigenfunction} then we say $T$ is a Type 2 tree. 

If $T$ is a Type 1 tree, then the only null vertex adjacent to a
non-null vertex (see Theorem \ref{thm:monoEigenfunction}) is said
to be the characteristic vertex of $T$.

If $T$ is a Type 2 tree, by Theorem \ref{thm:monoEigenfunction}
the only mixed block is formed by only two adjacent vertices. For
a Type 2 tree, we say that two vertices $i$ and $j$ are characteristic
vertices if and only if they are adjacent and satisfies $sign(g(i))\neq sign(g(j))$. 
\begin{thm}
\label{thm:type1Perron} Let $T$ be a tree and $g$ a eigenfunction
of $\lambda_{1}$. Then $T$ is a Type 1 tree with characteristic
vertex $v$ if and only if there are at least two Perron branches
at $v$. In this case, $\lambda_{1}=\dfrac{1}{\rho(\mathcal{L}(C)^{-1})}$
for each Perron branch $C$ at $v$. \end{thm}
\begin{proof}
Suppose that $T$ is a Type 1 tree and $v$ is its characteristic
vertex. Let $C_{0},C_{1},\ldots,C_{r}$ be the branchs of $T\setminus v$
and assume the normalized Laplacian matrix is in the form 
\begin{equation}
\mathcal{L}=\left[\begin{array}{ccccc}
\mathcal{L}(C_{0}) & 0 & \cdots & 0 & c_{0}\\
0 & \mathcal{L}(C_{1}) &  & 0 & c_{1}\\
\vdots &  & \ddots & \vdots & \vdots\\
0 & 0 & \cdots & \mathcal{L}(C_{r}) & c_{r}\\
\left(c_{0}\right)^{T} & \left(c_{1}\right)^{T} & \cdots & \left(c_{r}\right)^{T} & d_{v}
\end{array}\right],\label{matL}
\end{equation}
where $\mathcal{L}(C_{i})$ corresponds to vertices of the connected
component $C_{i}$, for $i=0,1,\ldots,r$ and $c_{i}$ is a 0,-1 vector
that accounts for the edges between the vertex $v$ and the connected
component $C_{i}$. For convenience, we assume that the last rows
and columns of $\mathcal{L}$ represent the vertex $v$.

We can define functions $g^{(i)}$ over each branch $C_{i}$ as $g^{(i)}(x)=g(x)$
where $x\in C_{i}$. From the relation $\mathcal{L}g=\lambda_{1}g$,
we have 
\[
\mathcal{L}(C_{i})g^{(i)}=\lambda_{1}y^{(i)},
\]
 since $g(v)=0$. From the fact that $\sum\sqrt{d_{x}}g(x)=0$, we
know $g$ assumes positive and negative values. Applying case 1 of
Theorem \ref{thm:monoEigenfunction}, we notice that there are functions
$g^{(r)}>0$ e $g^{(s)}<0$. Using Perron-Frobenius theorem, the only
positive eigenfunction are the Perron vector. Hence, $g^{(r)}$ and
$g^{(s)}$ are Perron vectors for $\mathcal{L}(C_{r})^{-1}$ and $\mathcal{L}(C_{s})^{-1}$,
respectively. If $g^{(i)}(x)=0$ for some $x\in C_{i}$, by case 1
of Theorem \ref{thm:monoEigenfunction} $g^{(i)}\equiv0$ and it is
not a Perron branch. Therefore, for each non-null branch we have the
relation for the Perron vector $g^{(i)}$ 
\[
\mathcal{L}(C_{i})^{-1}g^{(i)}=\frac{1}{\lambda_{1}(G)}g^{(i)}.
\]

It remains to show that $C_{r}$ and $C_{s}$ are Perron branch at
$v$. Suppose, by contradiction, that it is not true. Then it would
exist another component, say $C_{3}$, such that the Perron value
is larger than $1/\lambda_{1}$. We call $z$ the Perron vector of
$\mathcal{L}(C_{3})^{-1}$, normalized so that $\textbf{1}^{T}D_{3}^{\frac{1}{2}}z=1/\sqrt{2}$,
where $D_{3}$ is the diagonal degree matrix of $C_{3}$. Also, we
define $u=g^{(r)}/\sqrt{2}\textbf{1}^{T}D_{r}g^{(r)}$, where $D_{r}$
is the diagonal degree matrix of $C_{r}$. Thereof, we consider the
vector 
\[
w=[u,0,\cdots0,-z,0,\ldots,0]^{T},
\]
which is obviously orthogonal to $D\mathbf{1}$, where $D$ is the
diagonal degree matrix of $T$, and also $\Vert w\Vert=1$. Since 

\[
w^{T}\mathcal{L}w=\lambda_{1}u^{T}u+\frac{1}{\rho(\mathcal{L}(C_{3})^{-1})}z^{T}z<\lambda_{1}u^{T}u+\lambda_{1}z^{T}z=\lambda_{1}ww^{T}
\]
we obtain a contradiction with the fact that 
\[
\lambda_{1}(T)=\min_{\substack{\Vert x\Vert=1\\
x\bot D\textbf{1}
}
}x^{T}\mathcal{L}x.
\]
Thus, we obtain that $C_{r}$ and $C_{s}$ are indeed the Perron branchs
at $v$. This concludes the first part.

Conversely, assume that there are at least two Perron branchs at vertex
$v$, let us say $C_{i}$ and $C_{j}$ are two of them. Let $y$ and
$z$ be the Perron vectors of $\mathcal{L}(C_{i})^{-1}$ and $\mathcal{L}(C_{j})^{-1}$,
respectively. Taking into account (\ref{matL}), we can make $y$
and $z$ normalized such that $c_{i}^{T}x-c_{j}^{T}y=0$. Now, we
define the function $g$ as
\[
\begin{cases}
g(u)=y(u) & u\in C_{i};\\
g(u)=-z(u) & u\in C_{j};\\
0 & otherwise.
\end{cases}
\]
Hence, we have the relation $\mathcal{L}g=\frac{1}{\rho(\mathcal{L}(C_{i})^{-1})}g$.
It is easy to see that if $\lambda_{1}(T)=\frac{1}{\rho(\mathcal{L}(C_{i})^{-1})}$,
then $g$ is an eigenfunction that  makes $T$ a Type 1 tree with
characteristic vertex $v$, since $g(v)=0$. 

In order too see that $\lambda_{1}(T)=\frac{1}{\rho(\mathcal{L}(C_{i})^{-1})}$,
consider the submatrix $M$ of $\mathcal{L}$ obtained by deleting
the column and row corresponding to vertex $v$. It is easy to see
that the eigenvalues of $M$ are the union of eigenvalues of all matrices
$\mathcal{L}(C_{t})$, for $t=0,1,\ldots,r$. Since $\rho(\mathcal{L}(C_{i})^{-1})$
is the largest among all branchs, we can say that $\frac{1}{\rho(\mathcal{L}(C_{i})^{-1})}$
is the smallest eigenvalue of $M$. In fact, we have at least two
eigenvalues equal to $\frac{1}{\rho(\mathcal{L}(C_{i})^{-1})}$. Therefore,
we by the interlacing property of eigenvalues for principal submatrices,
we obtain 
\[
\frac{1}{\rho(\mathcal{L}(C_{i})^{-1})}\leq\lambda_{1}(T)\leq\frac{1}{\rho(\mathcal{L}(C_{i})^{-1})}.
\]
This shows the theorem.
\end{proof}
The previous theorem is a natural application of the same method used
in \cite{key-3} where, in the context of Laplacian matrix, it was
characterized the algebraic connectivity for Type I trees. However,
if we want to find some characterization for Type 2 trees using the
normalized Laplacian matrix, we must perform a different calculation
in order to obtain matrices that characterize $\lambda_{1}$. As the
next theorem show us, these matrices are more complicated than those
in \cite{key-3}.
\begin{thm}
\label{thm:type2Perron}Let $T$ be a tree on $n$ vertices with normalized
Laplacian matrix $\mathcal{L}$ and let $i$ and $j$ be adjacent
vertices of $T$. For $i$ and $j$ be characteristic vertices of
$T$ it is necessary and sufficient that there exists a $\gamma\in(0,1)$
such that 
\[
\rho(M_{1}-\gamma D_{1}^{\frac{1}{2}}\mathbf{1}\mathbf{1}^{T}D_{1}^{\frac{1}{2}})=\rho(M_{2}-\left(1-\gamma\right)D_{2}^{\frac{1}{2}}\mathbf{1}\mathbf{1}^{T}D_{2}^{\frac{1}{2}})=\frac{1}{\lambda_{1}},
\]
where $M_{1}$ is the normalized bottleneck matrix for the branch
at $j$ containing $i$ and $D_{1}$ is the degree matrix of this
branch; $M_{2}$ is the normalized bottleneck matrix for the branch
at $i$ containing $j$ and $D_{2}$ is the degree matrix of this
branch. \end{thm}
\begin{proof}
We can put the normalized Laplacian matrix of $T$ in the following
format
\[
\mathcal{L}=\left[\begin{array}{cc}
M_{1}^{-1} & -\frac{1}{\sqrt{d_{i}d_{j}}}e_{k}e_{1}^{T}\\
-\frac{1}{\sqrt{d_{i}d_{j}}}e_{1}e_{k}^{T} & M_{2}^{-1}
\end{array}\right],
\]
where the last row of $M_{1}^{-1}$ represents the vertex $i$ and
the first row of $M_{2}^{-1}$ represents the vertex $j$. 

First, we suppose that $i$ and $j$ are characteristic vertices of
$T$ . By using part (1) of the Theorem \ref{thm:valuationEigenfunction},
we have that both branchs at $i$ and at $j$ have the same sign each.
Moreover, the theorem ensures that we can write the eigenvector associated
with $\lambda_{1}$ as $v=\left[-v_{1}|v_{2}\right]^{T}$, where $v_{1}$and
$v_{2}$ are both positive vectors. Since $\mathbf{1}^{T}D^{\frac{1}{2}}v=0$,
we have $\mathbf{1}^{T}D_{1}^{\frac{1}{2}}v_{1}=\mathbf{1}^{T}D_{2}^{\frac{1}{2}}v_{1}$.

From the equation $\mathcal{L}v=\lambda_{1}v$, if we set $\alpha=e_{1}^{T}v_{2}$
and $\beta=e_{k}^{T}v_{1}$, we find that
\[
-M_{1}^{-1}v_{1}-\frac{\alpha}{\sqrt{d_{i}d_{j}}}e_{k}=-\lambda_{1}v_{1},
\]
which we can rewrite as

\[
\frac{v_{1}}{\lambda_{1}}=M_{1}v_{1}-\frac{\alpha}{\lambda_{1}\sqrt{d_{i}d_{j}}}M_{1}e_{k}.
\]
Using Lemma \ref{lem:bottleNorm}, we conclude that $M_{1}e_{k}=\sqrt{d_{i}}D_{1}^{\frac{1}{2}}\mathbf{1}$,
because $\left|P_{a,i,j}\right|=1$ for any vertex $a$ in the branch
at $j$ containing $i$ . Hence, we have
\begin{equation}
\frac{v_{1}}{\lambda_{1}}=M_{1}v_{1}-\frac{\alpha}{\lambda_{1}\sqrt{d_{j}}}D_{1}^{\frac{1}{2}}\mathbf{1}.\label{eq:Eq1}
\end{equation}
Now we multiply $e_{k}^{T}$ by (\ref{eq:Eq1}), to obtain 
\begin{eqnarray*}
\frac{e_{k}^{T}v_{1}}{\lambda_{1}} & = & e_{k}^{T}\left(M_{1}v_{1}-\frac{\alpha}{\lambda_{1}\sqrt{d_{j}}}D_{1}^{\frac{1}{2}}\mathbf{1}\right)\\
 & = & \sqrt{d_{i}}\mathbf{1}^{T}D_{1}^{\frac{1}{2}}v_{1}-\frac{\alpha\sqrt{d_{i}}}{\lambda_{1}\sqrt{d_{j}}}.
\end{eqnarray*}
Hence, we obtain 
\[
\frac{\beta}{\lambda_{1}}=\sqrt{d_{i}}\mathbf{1}^{T}D_{1}^{\frac{1}{2}}v_{1}-\frac{\alpha\sqrt{d_{i}}}{\lambda_{1}\sqrt{d_{j}}}
\]
which can be rewritten as 
\[
\frac{\sqrt{d_{i}d_{j}}}{\beta\sqrt{d_{j}}+\alpha\sqrt{d_{i}}}\mathbf{1}^{T}D_{1}^{\frac{1}{2}}v_{1}=\frac{1}{\lambda_{1}}.
\]
Now, we replace $\frac{1}{\lambda_{1}}$ in (\ref{eq:Eq1}), to obtain
\begin{eqnarray*}
\frac{v_{1}}{\lambda_{1}} & = & M_{1}v_{1}-\frac{\alpha\sqrt{d_{i}d_{j}}}{\sqrt{d_{j}}\left(\beta\sqrt{d_{j}}+\alpha\sqrt{d_{i}}\right)}\mathbf{1}^{T}D_{1}^{\frac{1}{2}}v_{1}D_{1}^{\frac{1}{2}}\mathbf{1}\\
 & = & M_{1}v_{1}-\frac{\alpha\sqrt{d_{i}}}{\beta\sqrt{d_{j}}+\alpha\sqrt{d_{i}}}D_{1}^{\frac{1}{2}}\mathbf{1}\mathbf{1}^{T}D_{1}^{\frac{1}{2}}v_{1}
\end{eqnarray*}
Therefore, we have 
\[
\frac{v_{1}}{\lambda_{1}}=\left(M_{1}-\frac{\alpha\sqrt{d_{i}}}{\beta\sqrt{d_{j}}+\alpha\sqrt{d_{i}}}D_{1}^{\frac{1}{2}}\mathbf{1}\mathbf{1}^{T}D_{1}^{\frac{1}{2}}\right)v_{1.}
\]

The same calculation for the matrix $M_{2}$, give us the relation
\[
\frac{v_{2}}{\lambda_{1}}=\left(M_{2}-\frac{\beta\sqrt{d_{j}}}{\beta\sqrt{d_{j}}+\alpha\sqrt{d_{i}}}D_{2}^{\frac{1}{2}}\mathbf{1}\mathbf{1}^{T}D_{2}^{\frac{1}{2}}\right)v_{1.}
\]

Now, from Lemma (\ref{lem:bottleNorm}), we conclude that $\left|P_{a,b,j}\right|\geq1$
and $\left|P_{a,b,i}\right|\geq1$, since the edge between $i$ and
$j$ is in any set of that form. Hence, we have $M_{1}\geq D_{1}^{\frac{1}{2}}\mathbf{1}\mathbf{1}^{T}D_{1}^{\frac{1}{2}}$
and $M_{2}\geq D_{2}^{\frac{1}{2}}\mathbf{1}\mathbf{1}^{T}D_{2}^{\frac{1}{2}}$.
Besides, if we define $\gamma=\frac{\beta\sqrt{d_{j}}}{\beta\sqrt{d_{j}}+\alpha\sqrt{d_{i}}}$
and notice that $\gamma\in(0,1)$, we conclude that $v_{1}$ is a
positive eigenvector of the positive matrix $M_{1}-\gamma D_{1}^{\frac{1}{2}}\mathbf{1}\mathbf{1}^{T}D_{1}^{\frac{1}{2}}$
and that $v_{2}$is a positive eigenvector for the matrix $M_{2}-\left(1-\gamma\right)D_{2}^{\frac{1}{2}}\mathbf{1}\mathbf{1}^{T}D_{2}^{\frac{1}{2}}$.
Therefore, from the Perron-Frobenius theory, we have 
\[
\rho(M_{1}-\gamma D_{1}^{\frac{1}{2}}\mathbf{1}\mathbf{1}^{T}D_{1}^{\frac{1}{2}})=\rho(M_{2}-\left(1-\gamma\right)D_{2}^{\frac{1}{2}}\mathbf{1}\mathbf{1}^{T}D_{2}^{\frac{1}{2}})=\frac{1}{\lambda_{1}},
\]
as required.

Reciprocally, assume that there is a $\gamma\in(0,1)$ that satisfies
$\rho(M_{1}-\gamma D_{1}^{\frac{1}{2}}\mathbf{1}\mathbf{1}^{T}D_{1}^{\frac{1}{2}})=\rho(M_{2}-\left(1-\gamma\right)D_{2}^{\frac{1}{2}}\mathbf{1}\mathbf{1}^{T}D_{2}^{\frac{1}{2}})=\frac{1}{\lambda_{1}},$
where $v_{1}$ and $v_{2}$ are the Perron vectors of $M_{1}-\gamma D_{1}^{\frac{1}{2}}\mathbf{1}\mathbf{1}^{T}D_{1}^{\frac{1}{2}}$
and $M_{2}-\left(1-\gamma\right)D_{2}^{\frac{1}{2}}\mathbf{1}\mathbf{1}^{T}D_{2}^{\frac{1}{2}}$,
respectively. Then we can compute
\begin{eqnarray*}
\frac{e_{k}^{T}v_{1}}{\lambda_{1}} & = & e_{k}^{T}\left(M_{1}-\gamma D_{1}^{\frac{1}{2}}\mathbf{1}\mathbf{1}^{T}D_{1}^{\frac{1}{2}}\right)v_{1}\\
 & = & \left(\sqrt{d_{i}}\mathbf{1}^{T}D_{1}^{\frac{1}{2}}-\gamma\sqrt{d_{i}}\mathbf{1}^{T}D_{1}^{\frac{1}{2}}\right)v_{1}\\
 & = & \left(1-\gamma\right)\sqrt{d_{i}}\mathbf{1}^{T}D_{1}^{\frac{1}{2}}v_{1}.
\end{eqnarray*}
Also, we can choose the eigenvectors $v_{1}$ and $v_{2}$ normalized
such that $\mathbf{1}^{T}D_{1}^{\frac{1}{2}}v_{1}=\mathbf{1}^{T}D_{2}^{\frac{1}{2}}v_{1}$,
and then we can write
\begin{eqnarray}
\frac{e_{k}^{T}v_{1}}{\lambda_{1}} & = & \left(1-\gamma\right)\sqrt{d_{i}}\mathbf{1}^{T}D_{2}^{\frac{1}{2}}v_{2}.\label{eq:eq2}
\end{eqnarray}

Similarly, using the same procedure, we can compute 
\[
e_{1}^{T}\left(M_{2}-\left(1-\gamma\right)D_{2}^{\frac{1}{2}}\mathbf{1}\mathbf{1}^{T}D_{2}^{\frac{1}{2}}\right)v_{2}
\]
 to obtain the relation 
\begin{equation}
\frac{e_{1}^{T}v_{2}}{\lambda_{1}}=\gamma\sqrt{d_{j}}\mathbf{1}^{T}D_{1}^{\frac{1}{2}}v_{1}.\label{eq:eq2b}
\end{equation}

Using the relation (\ref{eq:eq2}) in the equation $\left(M_{2}-\left(1-\gamma\right)D_{2}^{\frac{1}{2}}\mathbf{1}\mathbf{1}^{T}D_{2}^{\frac{1}{2}}\right)v_{2}=\frac{1}{\lambda_{1}}v_{2}$,
we obtain
\begin{eqnarray*}
\frac{1}{\lambda_{1}}v_{2} & = & M_{2}v_{2}-\left(1-\gamma\right)D_{2}^{\frac{1}{2}}\mathbf{1}\mathbf{1}^{T}D_{2}^{\frac{1}{2}}v_{2}\\
 & = & M_{2}v_{2}-\frac{1}{\lambda_{1}\sqrt{d_{i}}}D_{2}^{\frac{1}{2}}\mathbf{1}e_{k}^{T}v_{1}.
\end{eqnarray*}
By applying Lemma \ref{lem:bottleNorm}, we use the relation $M_{2}e_{1}=\sqrt{d_{j}}D_{2}^{\frac{1}{2}}\mathbf{1}$
to get
\[
\frac{1}{\lambda_{1}}v_{2}=M_{2}v_{2}-\frac{1}{\lambda_{1}\sqrt{d_{i}d_{j}}}M_{2}e_{1}e_{k}^{T}v_{1},
\]
which is equivalent to 
\begin{equation}
\lambda_{1}v=M_{2}^{-1}v_{2}+\frac{1}{\sqrt{d_{i}d_{j}}}e_{1}e_{k}^{T}v_{1}.\label{eq:eq3}
\end{equation}
In the same way, we can use the relation (\ref{eq:eq2b}) and then
rewrite the equation 
\[
\left(M_{1}-\gamma D_{1}^{\frac{1}{2}}\mathbf{1}\mathbf{1}^{T}D_{1}^{\frac{1}{2}}\right)v_{1}=\frac{1}{\lambda_{1}}v_{1}
\]
 as follows 
\begin{equation}
-\lambda_{1}v_{1}=-M_{1}^{-1}v_{1}+\frac{1}{\sqrt{d_{i}d_{j}}}e_{k}e_{1}^{T}v_{2}.\label{eq:eq4}
\end{equation}
Therefore, equation (\ref{eq:eq3}) and (\ref{eq:eq4}) show that
the vector $v=\left[-v_{1}|v_{2}\right]^{T}$ satisfies $\mathcal{L}v=\lambda_{1}v$.
This proofs the result.
\end{proof}

\section{Normalized Bottleneck Matrix}

The previous section pointed us to the bottleneck matrices in order
to characterize $\lambda_{1}$ of trees. Hence, in this section we
shall perform a more careful analysis on the structure of these matrices
with the expectation of giving prolific results about $\lambda_{1}$.
In fact, it allows us to extremize the $\lambda_{1}$ over the set
of trees. 

First, we define the set $P_{i,j,k}$ as the set of edges of $T$
which are on both the path from vertex $i$ to vertex $k$ and the
path from the vertex $j$ to vertex $k$. The following lemma was
obtained by Kirkland in \cite{key-3}, where it was investigated Perron
components of trees using the Laplacian matrix. 
\begin{lem}
\label{lem:bottleLaplacian}Consider a tree $T$ at $n$ vertex. Denote
by $L_{k}$ the principal submatrix of the Laplacian matrix $L\left(T\right)$
obtained by deleting the $k-$th column and the $k-$the row from
$L\left(T\right)$. Then the entry $\left(i,j\right)$ of $L_{k}^{-1}$
equals to the number of edges at $P_{i,j,k}$.
\end{lem}
The following lemma concerns the normalized Laplacian, and also we
can describe the entries of $\mathcal{L}_{k}^{-1}$.
\begin{lem}
\label{lem:bottleNorm}Consider a tree $T$ with $n$ vertex. Then
$\left(i,j\right)$ entry of $\mathcal{L}_{k}^{-1}$ is equal to $\sqrt{d_{i}d_{j}}\left|P_{i,j,k}\right|$.\end{lem}
\begin{proof}
We observe that, since $D$ is a diagonal matrix, then 
\[
\mathcal{L}_{k}=\left(D^{-\frac{1}{2}}LD^{-\frac{1}{2}}\right)_{k}=D_{k}^{-\frac{1}{2}}L_{k}D_{k}^{-\frac{1}{2}}.
\]

Thus, it is straightforward to obtain $\mathcal{L}_{k}^{-1}=D_{k}^{\frac{1}{2}}L_{k}^{-1}D_{k}^{\frac{1}{2}}$.
By applying Lemma \ref{lem:bottleLaplacian}, we obtain that the $\left(i,j\right)$
entry of $\mathcal{L}_{k}^{-1}$ is equal to 

\[
\left(D_{k}^{\frac{1}{2}}L_{k}^{-1}D_{k}^{\frac{1}{2}}\right)_{i,j}=\left(D_{k}^{\frac{1}{2}}\right)_{i,i}\left|P_{i,j,k}\right|\left(D_{k}^{\frac{1}{2}}\right)_{j,j}=\sqrt{d_{i}d_{j}}\left|P_{i,j,k}\right|.
\]

\end{proof}
The next result describes Perron branchs of trees in a similar fashion
to \cite{key-3}.
\begin{lem}
\label{lem:type2branch}$T$ is a Type 2 tree with characteristic
$i$ and $j$ if and only if $i$ and $j$ are adjacent and the branch
at $i$ containing vertex $j$ is the unique Perron branch at $i$,
while the branch at $j$ containing $i$ is the unique Perron branch
at $j$.\end{lem}
\begin{proof}
By Theorem \ref{thm:type2Perron}, for $T$ be a Type 2 it is necessary
and sufficient that there exists a $\gamma\in(0,1)$ such that 
\[
\rho(M_{1}-\gamma D_{1}^{\frac{1}{2}}\mathbf{1}\mathbf{1}^{T}D_{1}^{\frac{1}{2}})=\rho(M_{2}-\left(1-\gamma\right)D_{2}^{\frac{1}{2}}\mathbf{1}\mathbf{1}^{T}D_{2}^{\frac{1}{2}})=\frac{1}{\lambda_{1}}.
\]
We consider the values $\rho(M_{1}-xD\mathbf{1}\mathbf{1}^{T}D_{1}^{\frac{1}{2}})$
and $\rho(M_{2}-\left(1-x\right)D_{2}^{\frac{1}{2}}\mathbf{1}\mathbf{1}^{T}D_{2}^{\frac{1}{2}})$
as functions of $x$, and notice that they are decreasing and increasing,
respectively. Also, they are both continuous functions, hence the
existence of such $\gamma\in(0,1)$ is equivalent to $\rho(M_{1}-D\mathbf{1}\mathbf{1}^{T}D_{1}^{\frac{1}{2}})<\rho(M_{2})$
and $\rho(M_{1})>\rho(M_{2}-D_{2}^{\frac{1}{2}}\mathbf{1}\mathbf{1}^{T}D_{2}^{\frac{1}{2}})$. 

Using the description of the entries of $M_{1}$ and $M_{2}$ given
by Lemma \ref{lem:bottleNorm}, it is easy to see that the matrix
$M_{2}-D_{2}^{\frac{1}{2}}\mathbf{1}\mathbf{1}^{T}D_{2}^{\frac{1}{2}}$
is similar to the bottleneck matrix of the branchs at $j$ which do
not contain $i$. Also, we the matrix $M_{1}-D\mathbf{1}\mathbf{1}^{T}D_{1}^{\frac{1}{2}}$
is similar to the bottleneck matrix of the branchs at $i$ which do
not contain $j$. Therefore, the inequalities $\rho(M_{1}-D\mathbf{1}\mathbf{1}^{T}D_{1}^{\frac{1}{2}})<\rho(M_{2})$
and $\rho(M_{1})>\rho(M_{2}-D_{2}^{\frac{1}{2}}\mathbf{1}\mathbf{1}^{T}D_{2}^{\frac{1}{2}})$
holds if and only if the branch at $i$ containing vertex $j$ is
the unique Perron branch at $i$, while the branch at $j$ containing
$i$ is the unique Perron branch at $j$.
\end{proof}
The following result provides a simple way to characterize Type 1
and Type 2 trees.
\begin{thm}
Let $T$ be a tree. $T$ is a Type 1 tree if and only if there is
only one vertex such that there are at least two Perron branchs. $T$
is a Type 2 tree if and only if at each vertex there is a unique Perron
branch.\end{thm}
\begin{proof}
First, assume that there is only one vertex such that there are least
two Perron branchs. Then by Theorem \ref{thm:type1Perron}, $T$ is
a Type I tree. Conversely, assume that $T$ is a Type 1 tree with
characteristic vertex $v$. Take any branch at some vertex $u\neq v$.
Let $P$ be the branch at $u$ containing $v$ and $Q$ be any other
branch at $u$. Let $C$ be the component at $v$ that contains $u$.
In light of Lemma \ref{lem:bottleNorm}, we can see that $\mathcal{L}(Q)^{-1}\leq\mathcal{L}(C)^{-1}\leq\mathcal{L}(P)^{-1}$
with the strict inequality in at least one entry. Hence we conclude
that $\rho(\mathcal{L}(Q)^{-1})<\rho(\mathcal{L}(P)^{-1})$ and that
there is only one Perron component at $u$.

If $T$ is a Type 2 tree, then by Lemma \ref{lem:type2branch} there
are a pair of adjacent vertex $i$ and $j$ such there is a unique
Perron branch at each one. If we consider a vertex different from
$i$ and $j$, then we can use the same argument of the previous part
to conclude that there is only one Perron branch at this vertex. Finally,
assume that at each vertex there is a unique Perron branch. If $T$
is not a Type 2 tree, then we have a contradiction with Theorem \ref{thm:type1Perron}.
This completes the theorem.\end{proof}
\begin{cor}
Let $T$ be a tree and $u$ a vertex which not be its characteristic
vertex. Then the unique Perron branch at $u$ is the branch containing
the characteristic vertex or vertices of $T$.\end{cor}

\end{document}